\documentclass[onecolumn]{autart}
\usepackage{amsmath}
\begin{document}
\input{amssym}
\begin{frontmatter}
\title{Symmetry classification of Newtonian incompressible fluid's equations flow in turbulent boundary layers}
\author[MN]{Mehdi Nadjafikhah}\ead{m\_nadjafikhah@iust.ac.ir},
\author[SRH]{Seyed Reza Hejazi}\ead{reza\_hejazi@iust.ac.ir}
\address[MN]{School of Mathematics, Iran University of Science and Technology, Narmak-16, Teharn, I.R.Iran.}
\address[SRH]{Same as author one.}
\begin{abstract}
Lie symmetry group method is applied to study Newtonian
incompressible fluid's equations flow in turbulent boundary
layers. The symmetry group and its optimal system are given, and
group invariant solutions associated to the symmetries are
obtained. Finally the structure of the Lie algebra symmetries is
determined.
\end{abstract}
\end{frontmatter}
\section{Introduction}
In physics and fluid mechanics, a  \textit{boundary layer} is that
layer of fluid in the immediate vicinity of a bounding surface. In
the Earth's atmosphere, the planetary boundary layer is the air
layer near the ground affected by diurnal heat, moisture or
momentum transfer to or from the surface. On an aircraft wing the
boundary layer is the part of the flow close to the wing. The
boundary layer effect occurs at the field region in which all
changes occur in the flow pattern. The boundary layer distorts
surrounding non-viscous flow. It is a phenomenon of viscous
forces. This effect is related to the Reynolds number (In fluid
mechanics and heat transfer, the \textit{Reynold's number} is a
dimensionless number that gives a measure of the ratio of inertial
forces to viscous and, consequently, it quantifies the relative
importance of these two types of forces for given flow
conditions). Laminar boundary layers come in various forms and can
be loosely classified according to their structure and the
circumstances under which they are created. The thin shear layer
which develops on an oscillating body is an example of a Stokes
boundary layer, whilst the Blasius boundary layer refers to the
well-known similarity solution for the steady boundary layer
attached to a flat plate held in an oncoming unidirectional flow.
When a fluid rotates, viscous forces may be balanced by the
Coriolis effect, rather than convective inertia, leading to the
formation of an \textit{Ekman} layer. Thermal boundary layers also
exist in heat transfer. Multiple types of boundary layers can
coexist near a surface simultaneously. The deduction of the
boundary layer equations was perhaps one of the most important
advances in fluid dynamics. Using an order of magnitude analysis,
the well-known governing \textsl{Navier$-$Stokes} equations of
viscous fluid flow can be greatly simplified within the boundary
layer. Notably, the characteristic of the partial differential
equations (PDE) becomes parabolic, rather than the elliptical form
of the full Navier$-$Stokes equations. This greatly simplifies the
solution of the equations. By making the boundary layer
approximation, the flow is divided into an inviscid portion (which
is easy to solve by a number of methods) and the boundary layer,
which is governed by an easier to solve PDE.

$~~~$Flow and heat transfer of an incompressible viscous fluid
over a stretching sheet appear in several manufacturing processes
of industry such as the extrusion of polymers, the cooling of
metallic plates, the aerodynamic extrusion of plastic sheets, etc.
In the glass industry, blowing, floating or spinning of fibres are
processes, which involve the flow due to a stretching surface.
Mahapatra and Gupta studied the steady two-dimensional
stagnation-point flow of an incompressible viscous fluid over a
flat deformable sheet when the sheet is stretched in its own plane
with a velocity proportional to the distance from the
stagnation-point. They concluded that, for a fluid of small
kinematic viscosity, a boundary layer is formed when the
stretching velocity is less than the free stream velocity and an
inverted boundary layer is formed when the stretching velocity
exceeds the free stream velocity. Temperature distribution in the
boundary layer is determined when the surface is held at constant
temperature giving the so called surface heat flux. In their
analysis, they used the finite-differences scheme along with the
\textit{Thomas algorithm} to solve the resulting system of
ordinary differential equations.

$~~~$The treatment of turbulent boundary layers is far more
difficult due to the time-dependent variation of the flow
properties. One of the most widely used techniques in which
turbulent flows are tackled is to apply \textit{Reynolds
decomposition}. Here the instantaneous flow properties are
decomposed into a mean and fluctuating component. Applying this
technique to the boundary layer equations gives the full turbulent
boundary layer equations
\begin{eqnarray}\label{eq:1}
\left\{ \begin{array}{l} \displaystyle {\frac{\partial
\bar{u}}{\partial x}+\frac{\partial \bar{v}}{\partial
y}=0},\\
\displaystyle {\bar{u}\frac{\partial \bar{u}}{\partial
x}+\bar{v}\frac{\partial \bar{u}}{\partial
y}=-\frac{1}{\rho}\frac{\partial \bar{p}}{\partial
x}+\nu\Big(\frac{\partial^2\bar{u}}{\partial
x^2}+\frac{\partial^2\bar{u}}{\partial
y^2}\Big)-\frac{\partial}{\partial
y}(\overline{u'v'})-\frac{\partial}{\partial
x}(\overline{u'^2})},\\
\displaystyle {\bar{u}\frac{\partial \bar{v}}{\partial
x}+\bar{v}\frac{\partial \bar{v}}{\partial
y}=-\frac{1}{\rho}\frac{\partial \bar{p}}{\partial
y}+\nu\Big(\frac{\partial^2\bar{v}}{\partial
x^2}+\frac{\partial^2\bar{v}}{\partial
y^2}\Big)-\frac{\partial}{\partial
x}(\overline{u'v'})-\frac{\partial}{\partial y}(\overline{v'^2})},
\end{array} \right.
\end{eqnarray}
where $\rho$ is the \textit{density}, $p$ is the
\textit{pressure}, $\nu$ is the \textit{kinematic viscosity} of
the fluid at a point and $\bar{u}$ and $\bar{v}$ are average of
the velocity components in Reynold decomposition. Here $u'$ and
$v'$ are the \textit{velocity fluctations} such that;
$u=\bar{u}+u'$ and $v=\bar{v}+v'$. By using the \textit{scale
analysis}(a powerful tool used in the mathematical sciences for
the simplification of equations with many terms), it can be shown
that the system (\ref{eq:1}) reduce to the classical form
\begin{eqnarray}\label{eq:22}
\left\{ \begin{array}{l} \displaystyle {\frac{\partial
\bar{u}}{\partial x}+\frac{\partial \bar{v}}{\partial
y}=0},\\
\displaystyle{\bar{u}\frac{\partial\bar{u}}{\partial
x}+\bar{v}\frac{\partial\bar{v}}{\partial
y}=-\frac{1}{\rho}\frac{\partial\bar{p}}{\partial
x}+\nu\frac{\partial^2\bar{u}}{\partial
y^2}-\frac{\partial}{\partial y}(\overline{u'v'})},\\
\displaystyle{\frac{\partial\bar{p}}{\partial y}=0}.
\end{array} \right.
\end{eqnarray}
The term $\overline{u'v'}$ in the system (\ref{eq:22}) called
\textit{Reynolds shear stress}, a tensor that conventionally
written
\begin{eqnarray}\label{eq:23}
R_{ij}\equiv\rho\overline{u'_iu'_j}.
\end{eqnarray}
The divergence of this stress is the force density on the fluid
due to the \textit{turbulent fluctuations}. Using Navier–Stokes
equations for a fluid whose stress versus rate of strain curve is
linear and passes through the origin (Newtonian fluid) the tensor
(\ref{eq:23}) reduces to
\begin{eqnarray*}
R_{ij}\equiv\mu\frac{\partial\bar{u}_i}{\partial x_j},
\end{eqnarray*}
where $\mu$ is the \textit{fluid viscosity}.

$~~~$This paper is concerned with the symmetry of Newtonian
incompressible fluid's equations flow in turbulent boundary layers
of the form (\ref{eq:22}). Lie-group theory is applied to the
equations of motion for determining symmetry reductions of partial
differential equations. The solution of the turbulent boundary
layer equations therefore necessitates the use of a turbulence
model, which aims to express the Reynolds shear stress in terms of
known flow variables or derivatives. The lack of accuracy and
generality of such models is a major obstacle in the successful
prediction of turbulent flow properties in modern fluid dynamics
\section{Lie Symmetries of the Equations}
A PDE with $p-$independent and $q-$dependent variables has a Lie
point transformations
\begin{eqnarray*}
\widetilde{x}_i=x_i+\varepsilon\xi_i(x,u)+{\mathcal
O}(\varepsilon^2),\qquad
\widetilde{u}_{\alpha}=u_\alpha+\varepsilon\varphi_\alpha(x,u)+{\mathcal
O}(\varepsilon^2),
\end{eqnarray*}
where
$\displaystyle{\xi_i=\frac{\partial\widetilde{x}_i}{\partial\varepsilon}\Big|_{\varepsilon=0}}$
for $i=1,...,p$ and
$\displaystyle{\varphi_\alpha=\frac{\partial\widetilde{u}_\alpha}{\partial\varepsilon}\Big|_{\varepsilon=0}}$
for $\alpha=1,...,q$. The action of the Lie group can be
considered by its associated infinitesimal generator
\begin{eqnarray}\label{eq:18}
\textbf{v}=\sum_{i=1}^p\xi_i(x,u)\frac{\partial}{\partial{x_i}}+\sum_{\alpha=1}^q\varphi_\alpha(x,u)\frac{\partial}{\partial{u_\alpha}}
\end{eqnarray}
on the total space of PDE (the space containing independent and
dependent variables). Furthermore, the characteristic of the
vector field (\ref{eq:18}) is given by
\begin{eqnarray*}
Q^\alpha(x,u^{(1)})=\varphi_\alpha(x,u)-\sum_{i=1}^p\xi_i(x,u)\frac{\partial
u^\alpha}{\partial x_i},
\end{eqnarray*}
and its $n-$th prolongation is determined by
\begin{eqnarray*}
\textbf{v}^{(n)}=\sum_{i=1}^p\xi_i(x,u)\frac{\partial}{\partial
x_i}+\sum_{\alpha=1}^q\sum_{\sharp
J=j=0}^n\varphi^J_\alpha(x,u^{(j)})\frac{\partial}{\partial
u^\alpha_J},
\end{eqnarray*}
where
$\varphi^J_\alpha=D_JQ^\alpha+\sum_{i=1}^p\xi_iu^\alpha_{J,i}$.
($D_J$ is the total derivative operator describes in
(\ref{eq:3})).

$~~~$Let us consider a one-parameter Lie group of infinitesimal
transformations $(x,y,\bar{u},\bar{v},\bar{p})$ given by
\begin{align*}
\widetilde{x}&=x+\varepsilon\xi_1(x,y,\bar{u},\bar{v},\bar{p})+{\mathcal
O}(\varepsilon^2),&
&&\widetilde{y}&=y+\varepsilon\xi_2(x,y,\bar{u},\bar{v},\bar{p})+{\mathcal
O}(\varepsilon^2),\\
\widetilde{\bar{u}}&=\bar{u}+\varepsilon\eta_1(x,y,\bar{u},\bar{v},\bar{p})+{\mathcal
O}(\varepsilon^2),&
&&\widetilde{\bar{v}}&=\bar{v}+\varepsilon\eta_2(x,y,\bar{u},\bar{v},\bar{p})+{\mathcal
O}(\varepsilon^2),\\
\widetilde{\bar{p}}&=\bar{p}+\varepsilon\eta_3(x,y,\bar{u},\bar{v},\bar{p})+{\mathcal
O}(\varepsilon^2),
\end{align*}
where $\varepsilon$ is the group parameter. Then one requires that
this transformations leaves invariant the set of solutions of the
system (\ref{eq:22}). This yields to the linear system of
equations for the infinitesimals
$\xi_1(x,y,\bar{u},\bar{v},\bar{p}),\xi_2(x,y,\bar{u},\bar{v},\bar{p})$,
$\eta_1(x,y,\bar{u},\bar{v},\bar{p}),\eta_2(x,y,\bar{u},\bar{v},\bar{p})$
and $\eta_3(x,y,\bar{u},\bar{v},\bar{p})$. The Lie algebra of
infinitesimal symmetries is the set of vector fields in the form
of
\begin{eqnarray*}
\textbf{v}=\xi_1\frac{\partial}{\partial
x}+\xi_2\frac{\partial}{\partial y}+\eta_1\frac{\partial}{\partial
\bar{u}}+\eta_2\frac{\partial}{\partial
\bar{v}}+\eta_3\frac{\partial}{\partial \bar{p}}.
\end{eqnarray*}
This vector field has the second prolongation
\begin{eqnarray}\label{eq:21}
\textbf{v}^{(2)}&=&\textbf{v}+\varphi_1^x\frac{\partial}{\partial
x}+\varphi_1^y\frac{\partial}{\partial
y}+\varphi_1^{xx}\frac{\partial}{\partial
\bar{u}_{xx}}+\varphi_1^{xy}\frac{\partial}{\partial
\bar{u}_{xy}}+\varphi_1^{yy}\frac{\partial}{\partial \bar{u}_{yy}}
+\varphi_2^x\frac{\partial}{\partial
\bar{v}_x}+\varphi_2^y\frac{\partial}{\partial
\bar{v}_y}+\varphi_2^{xx}\frac{\partial}{\partial
\bar{v}_{xx}}\nonumber\\
&&+\varphi_2^{xy}\frac{\partial}{\partial
\bar{v}_{xy}}+\varphi_2^{yy}\frac{\partial}{\partial
\bar{v}_{yy}}+\varphi_3^x\frac{\partial}{\partial
\bar{p}_x}+\varphi_3^y\frac{\partial}{\partial
\bar{p}_y}+\varphi_3^{xx}\frac{\partial}{\partial
\bar{p}_{xx}}+\varphi_3^{xy}\frac{\partial}{\partial
\bar{p}_{xy}}+\varphi_3^{yy}\frac{\partial}{\partial
\bar{p}_{yy}},\nonumber
\end{eqnarray}
with the coefficients
\begin{align*}
\varphi_1^x&=D_x(\varphi_1-\xi_1\bar{u}_x-\xi_2\bar{u}_y)+\xi_1\bar{u}_{xx}+\xi_2\bar{u}_{xy},&
\varphi_1^y&=D_y(\varphi_1-\xi_1\bar{u}_x-\xi_2\bar{u}_y)+\xi_1\bar{u}_{xy}+\xi_2\bar{u}_{yy},\\
\varphi_1^{xx}&=D_x^2(\varphi_1-\xi_1\bar{u}_x-\xi_2\bar{u}_y)+\xi_1\bar{u}_{xxx}+\xi_2\bar{u}_{xxy},&
\varphi_1^{xy}&=D_xD_y(\varphi_1-\xi_1\bar{u}_x-\xi_2\bar{u}_y)+\xi_1\bar{u}_{xxy}+\xi_2\bar{u}_{xxy},\\
\varphi_1^{yy}&=D_y^2(\varphi_1-\xi_1\bar{u}_x-\xi_2\bar{u}_y)+\xi_1\bar{u}_{xxy}+\xi_2\bar{u}_{xyy},&
\varphi_2^x&=D_x(\varphi_2-\xi_1\bar{v}_x-\xi_2\bar{v}_y)+\xi_1\bar{v}_{xx}+\xi_2\bar{v}_{xy},\\
\varphi_2^y&=D_y(\varphi_2-\xi_1\bar{v}_x-\xi_2\bar{v}_y)+\xi_1\bar{v}_{xy}+\xi_2\bar{v}_{yy},&
\varphi_2^{xx}&=D_x^2(\varphi_2-\xi_1\bar{v}_x-\xi_2\bar{v}_y)+\xi_1\bar{v}_{xxx}+\xi_2\bar{v}_{xxy},\\
\varphi_2^{xy}&=D_xD_y(\varphi_2-\xi_1\bar{v}_x-\xi_2\bar{v}_y)+\xi_1\bar{v}_{xxy}+\xi_2\bar{v}_{xxy},&
\varphi_2^{yy}&=D_y^2(\varphi_2-\xi_1\bar{v}_x-\xi_2\bar{v}_y)+\xi_1\bar{v}_{xxy}+\xi_2\bar{v}_{xyy},\\
\varphi_3^x&=D_x(\varphi_3-\xi_1\bar{p}_x-\xi_2\bar{p}_y)+\xi_1\bar{p}_{xx}+\xi_2\bar{p}_{xy},&
\varphi_3^y&=D_y(\varphi_3-\xi_1\bar{p}_x-\xi_2\bar{p}_y)+\xi_1\bar{p}_{xy}+\xi_2\bar{p}_{yy},\\
\varphi_3^{xx}&=D_x^2(\varphi_3-\xi_1\bar{p}_x-\xi_2\bar{p}_y)+\xi_1\bar{p}_{xxx}+\xi_2\bar{p}_{xxy},&
\varphi_3^{xy}&=D_xD_y(\varphi_3-\xi_1\bar{p}_x-\xi_2\bar{p}_y)+\xi_1\bar{p}_{xxy}+\xi_2\bar{p}_{xxy},\\
\varphi_3^{yy}&=D_y^2(\varphi_3-\xi_1\bar{p}_x-\xi_2\bar{p}_y)+\xi_1\bar{p}_{xxy}+\xi_2\bar{p}_{xyy},
\end{align*}
where the operators $D_x$ and $D_y$ denote the total derivative
with respect to $x$ and $y$:
\begin{eqnarray}\label{eq:3}
D_x=\frac{\partial}{\partial x}+\bar{u}_x\frac{\partial}{\partial
\bar{u}}+\bar{v}_x\frac{\partial}{\partial
\bar{v}}+\bar{p}_x\frac{\partial}{\partial \bar{p}}\cdots,\quad
D_y=\frac{\partial}{\partial y}+\bar{u}_y\frac{\partial}{\partial
\bar{u}}+\bar{v}_y\frac{\partial}{\partial
\bar{v}}+\bar{p}_y\frac{\partial}{\partial \bar{p}}\cdots.
\end{eqnarray}
Using the invariance condition, i.e., applying the second
prolongation $\textbf{v}^{(2)}$, to system (\ref{eq:22}), and by
solving the linear system
\begin{eqnarray*} \left\{\begin{array}{l} \displaystyle
\textbf{v}^{(2)}\Big(\frac{\partial \bar{u}}{\partial
x}+\frac{\partial \bar{v}}{\partial
y}\Big)=0,\\[3mm] \displaystyle \textbf{v}^{(2)}\Big(\bar{u}\frac{\partial\bar{u}}{\partial
x}+\bar{v}\frac{\partial\bar{v}}{\partial
y}+\frac{1}{\rho}\frac{\partial\bar{p}}{\partial
x}-\nu\frac{\partial^2\bar{u}}{\partial
y^2}+\frac{\partial}{\partial y}(\overline{u'v'})\Big)=0,\\[3mm]
\displaystyle \textbf{v}^{(2)}\Big(\frac{\partial\bar{p}}{\partial
y}\Big)=0,\end{array}\right.\quad \pmod{(\ref{eq:22})}
\end{eqnarray*}
the following system of 17 determining equations yields:
\begin{eqnarray*}
\begin{array}{lclclcl}
{\xi_1}_y=0,&&{\xi_1}_{\bar{u}}=0,&&{\xi_1}_{\bar{v}}=0,&&{\xi_1}_{\bar{v}}=0,\\
{\xi_1}_{xx}=0,&&{\xi_2}_x=0,&&{\xi_2}_{\bar{u}}=0,&&{\xi_2}_{\bar{v}}=0,\\
{\xi_2}_{\bar{p}}=0,&&{\xi_2}_{yy}=0,&&{\eta_1}_x=0,&&{\eta_1}_y=0,\\
{\eta_1}_{\bar{u}}=0,&&{\eta_1}_{\bar{v}}=0,&&{\eta_1}_{\bar{p}}=2{\xi_1}_x-4{\xi_2}_y,&&{\eta_2}={\bar{u}}({\xi_1}_x-2{\xi_2}_y),\\
{\eta_3}=-{\xi_2}_y.
\end{array}
\end{eqnarray*}
The solution of the above system gives the following coefficients
of the vector field $\textbf{v}$:
\begin{eqnarray*}
\begin{array}{lclclclcl}
\xi_1=C_1+C_4x,&&\xi_2=C_2+C_5y,&&\eta_1=C_4\bar{u}-2C_5\bar{u},&&\eta_2=-C_5\bar{v},&&\eta_3=C_3+2C_4\bar{p}-4C_5\bar{p},
\end{array}
\end{eqnarray*}
where $C_i$ for $i=1,...,5$ are arbitrary constants, thus the Lie
algebra ${\goth g}$ of the Newtonian incompressible fluid's
equations flow in turbulent boundary layers is spanned by the five
vector fields
\begin{eqnarray*}
\textbf{v}_1&=&\frac{\partial}{\partial x},\\
\textbf{v}_2&=&\frac{\partial}{\partial y},\\
\textbf{v}_3&=&\frac{\partial}{\partial \bar{p}},\\
\textbf{v}_4&=&x\frac{\partial}{\partial
x}+\bar{u}\frac{\partial}{\partial\bar{u}}+2\bar{p}\frac{\partial}{\partial\bar{p}},\\
\textbf{v}_5&=&y\frac{\partial}{\partial
y}-2\bar{u}\frac{\partial}{\partial\bar{u}}-\bar{v}\frac{\partial}{\partial
\bar{v}}-4\bar{p}\frac{\partial}{\partial\bar{p}},
\end{eqnarray*}
The commutation relations between these vector fields is given by
the table (\ref{table:1}), where entry in row $i$ and column $j$
representing $[\textbf{v}_i,\textbf{v}_j]$.
\begin{table}
\caption{Commutation relations of $\goth g$ }\label{table:1}
\vspace{-0.3cm}\begin{eqnarray*}\hspace{-0.75cm}\begin{array}{lclclclclc}
\hline
  [\,,\,]       &\hspace{1.1cm}\textbf{v}_1         &\hspace{0.5cm}\textbf{v}_2         &\hspace{0.5cm}\textbf{v}_3          &\hspace{0.5cm}\textbf{v}_4        &\hspace{0.5cm}\textbf{v}_5\\ \hline
  \textbf{v}_1  &\hspace{1.1cm} 0                   &\hspace{0.5cm} 0                   &\hspace{0.5cm}0                     &\hspace{0.5cm}\textbf{v}_1        &\hspace{0.5cm}0  \\
  \textbf{v}_2  &\hspace{1.1cm} 0                   &\hspace{0.5cm} 0                   &\hspace{0.5cm}0                     &\hspace{0.5cm}0                   &\hspace{0.5cm}\textbf{v}_2\\
  \textbf{v}_3  &\hspace{1.1cm} 0                   &\hspace{0.5cm} 0                   &\hspace{0.5cm}0                     &\hspace{0.5cm}2\textbf{v}_3       &\hspace{0.5cm}-4\textbf{v}_3\\
  \textbf{v}_4  &\hspace{1.1cm} -\textbf{v}_1       &\hspace{0.5cm} 0                   &\hspace{0.5cm}-2\textbf{v}_3        &\hspace{0.5cm}0                   &\hspace{0.5cm}0\\
  \textbf{v}_5  &\hspace{1.1cm} 0                   &\hspace{0.5cm} -\textbf{v}_2       &\hspace{0.5cm}4\textbf{v}_3         &\hspace{0.5cm}0                   &\hspace{0.5cm}0\\
  \hline\end{array}\end{eqnarray*}\end{table}

The one-parameter groups $G_i$ generated by the base of $\goth g$
are given in the following table.
\begin{eqnarray*}
\begin{array}{lcl}
G_1:(x+\varepsilon,y,\bar{u},\bar{v},\bar{p}),&& G_2:(x,y+\varepsilon,\bar{u},\bar{v},\bar{p}),\\
G_3:(x,y,\bar{u},\bar{v},\bar{p}+\varepsilon),&& G_4:(xe^\varepsilon,y,\bar{u}e^\varepsilon,\bar{v},\bar{p}e^{2\varepsilon})\\
G_5:(x,ye^\varepsilon,\bar{u}e^{-2\varepsilon},\bar{v}e^{-\varepsilon},\bar{p}e^{-4\varepsilon})
\end{array}
\end{eqnarray*}
Since each group $G_i$ is a symmetry group and if
$\bar{u}=f(x,y),\bar{v}=g(x,y)$ and $\bar{p}=h(x,y)$ are solutions
of the system (\ref{eq:22}), so are the functions
\begin{align*}
\bar{u}_1&=f(x+\varepsilon,y),&\bar{u}_1&=f(x,y+\varepsilon),&\bar{u}_1&=f(x,y),&\bar{u}_1&=e^{-\varepsilon}f(xe^\varepsilon,y),&\bar{u}=e^{2\varepsilon}f(x,ye^\varepsilon),\\
\bar{v}_1&=g(x+\varepsilon,y),&\bar{v}_1&=g(x,y+\varepsilon),&\bar{v}_1&=g(x,y),&\bar{v}_1&=g(xe^\varepsilon,y),&\bar{v}=e^\varepsilon g(x,ye^\varepsilon),\\
\bar{p}_1&=h(x+\varepsilon,y),&\bar{p}_1&=h(x,y+\varepsilon),&\bar{p}_1&=h(x,y)+\varepsilon,&\bar{p}_1&=e^{-2\varepsilon}h(xe^\varepsilon,y),&\bar{p}=e^{4\varepsilon}h(x,ye^\varepsilon),
\end{align*}
where $\varepsilon$ is a real number. Here we can find the general
group of the symmetries by considering a general linear
combination $c_1\textbf{v}_1+\cdots+c_1\textbf{v}_5$ of the given
vector fields. In particular if $g$ is the action of the symmetry
group near the identity, it can be represented in the form
$g=\exp(\varepsilon\textbf{v}_5)\circ\cdots\circ\exp(\varepsilon\textbf{v}_1)$.
Consequently, if $\bar{u}=f(x,y),\bar{v}=g(x,y)$ and
$\bar{p}=h(x,y)$ be a solution for system (\ref{eq:22}), so is
\begin{eqnarray}\label{eq:13}
\bar{u}&=&e^\varepsilon
f\Big((x+\varepsilon)e^\varepsilon,(y+\varepsilon)e^\varepsilon\Big),\nonumber\\
\bar{v}&=&g\Big((x+\varepsilon)e^\varepsilon,(y+\varepsilon)e^\varepsilon\Big),\\
\bar{p}&=&e^{2\varepsilon}h\Big((x+\varepsilon)e^\varepsilon,(y+\varepsilon)e^\varepsilon\Big)+\varepsilon
e^{-\varepsilon},\nonumber
\end{eqnarray}
a fact that can, of course, be checked directly.
\section{Symmetry reduction for Newtonian incompressible fluid's equations flow in turbulent boundary layers}
Lie-group method is applicable to both linear and non-linear
partial differential equations, which leads to similarity
variables that may be used to reduce the number of independent
variables in partial differential equations. By determining the
transformation group under which a given partial differential
equation is invariant, we can obtain information about the
invariants and symmetries of that equation.

$~~~$The first advantage of symmetry group method is to construct
new solutions from known solutions. Neither the first advantage
nor the second will be investigated here, but symmetry group
method will be applied to the Eq. (\ref{eq:22}) to be connected
directly to some order differential equations. To do this, a
particular linear combinations of infinitesimals are considered
and their corresponding invariants are determined. The Newtonian
incompressible fluid's equations flow in turbulent boundary layers
expressed in the coordinates $(x,y)$, so to reduce this equation
is to search for its form in specific coordinates. Those
coordinates will be constructed by searching for independent
invariants $(r,s)$ corresponding to an infinitesimal generator. So
using the chain rule, the expression of the equation in the new
coordinate allows us to the reduced equation. Here we will obtain
some invariant solutions with respect to symmetries in the way of
similarity transformations to leave invariant the system
(\ref{eq:22}) under the group of symmetries. First we obtain the
similarity variables for each term of the Lie algebra $\goth g$,
then we use this method to reduced the PDE and find the invariant
solutions. Here our system has two-independent and four-dependent
variables, thus, the similarity transformations and invariant
transformations including of invariants up to order 2 comining in
table (\ref{table:2}) and (\ref{table:3}),
\begin{table}
\caption{Reduction of system}\label{table:2}
\vspace{-0.3cm}\begin{eqnarray*}\hspace{-0.75cm}\begin{array}{lclclclclclc}
\hline
  \mbox{vector field}  &\hspace{-1cm}\mbox{symmetry transformations}                                                                     &\hspace{2cm}\mbox{similarities}         \\ \hline
  \textbf{v}_1         &\hspace{-0.5cm}r=x+\varepsilon,s=y,f(r,s)=u,g(r,s)=v,h(r,s)=p                                                    &\hspace{0.5cm} x=s,y=r,u=f(r),v=g(r),p=h(r) \\
  \textbf{v}_2         &\hspace{-0.5cm}r=x,s=y+\varepsilon,f(r,s)=u,g(r,s)=v,h(r,s)=p                                                    &\hspace{0.5cm} x=r,y=s,u=f(r),v=g(r),p=h(r)\\
  \textbf{v}_3         &\hspace{-0.5cm}r=x,s=y,f(r,s)=u,g(r,s)=v,h(r,s)=p+\varepsilon                                                    &\hspace{0.5cm} \mbox{translation on $p$ is the similarity } \\
  \textbf{v}_4         &\hspace{-0.5cm}r=xe^\varepsilon,s=y,f(r,s)=e^{-\varepsilon}u,g(r,s)=v,h(r,s)=e^{-2\varepsilon}p                  &\hspace{0.5cm}x=e^s,y=r,u=e^sf(r),v=g(r),p=e^{2s}h(r) \\
  \textbf{v}_5         &\hspace{-0.5cm}r=x,s=ye^\varepsilon,f(r,s)=e^{2\varepsilon}u,g(r,s)=e^{\varepsilon}v,h(r,s)=e^{4\varepsilon}p    &\hspace{0.5cm}x=r,y=e^s,u=e^{-2s}f(r),v=e^{-s}g(r),p=e^{-4s}h(r) \\
  \hline\end{array}\end{eqnarray*}\end{table}
  \begin{table}
\caption{Invariants}\label{table:3}
\vspace{-0.3cm}\begin{eqnarray*}\hspace{-0.75cm}\begin{array}{lclclclclclc}
\hline
  \mbox{vector field}  &\hspace{0.5cm}\mbox{ordinary invariant}              &\hspace{1cm}\mbox{1st order diff. invariant}                              &\hspace{1cm}\mbox{2nd order diff. invariant}      \\ \hline
  \textbf{v}_1         &\hspace{-0.5cm}y,u,v,p                               &\hspace{0.5cm}*,u_x,v_x,p_x,u_y,v_y,p_y                                   &\hspace{0.5cm}**,u_{xx},v_{xx},p_{xx},u_{xy},v_{xy},p_{xy},u_{yy},v_{yy},p_{yy}\\
  \textbf{v}_2         &\hspace{-0.5cm}x,u,v,p                               &\hspace{0.5cm}*,u_x,v_x,p_x,u_y,v_y,p_y                                   &\hspace{0.5cm}**,u_{xx},v_{xx},p_{xx},u_{xy},v_{xy},p_{xy},u_{yy},v_{yy},p_{yy}\\
  \textbf{v}_3         &\hspace{-0.5cm}x,y,u,v                               &\hspace{0.5cm}*,u_x,v_x,p_x,u_y,v_y,p_y                                   &\hspace{0.5cm}**,u_{xx},v_{xx},p_{xx},u_{xy},v_{xy},p_{xy},u_{yy},v_{yy},p_{yy}\\
  \textbf{v}_4         &\hspace{-0.5cm}u,\frac{u}{x},v,\frac{p}{x^2}         &\hspace{0.5cm}*,u_x,xv_x,\frac{p_x}{x},\frac{u_y}{x},v_y,\frac{p_y}{x^2}  &\hspace{0.5cm}y,\frac{u}{x},v,p,u_x,xv_x,\frac{p_x}{x},\frac{u_y}{x},v_y,\frac{p_y}{x^2},\\
                                                                                                                                                        &&&xu_{xx},x^2v_{xx},p_{xx},u_{xy},xv_{xy},\frac{p_{xy}}{x},\frac{u_{yy}}{x},v_{yy},\frac{p_{yy}}{x^2} \\
  \textbf{v}_5         &\hspace{-0.5cm}x,y^2u,yv,y^4p                        &\hspace{0.5cm}*,y^2u_x,yv_x,y^4p_x,y^3u_y,y^2v_y,y^5p_y                   &\hspace{0.5cm}*,y^2u_x,yv_x,y^4p_x,y^3u_y,y^2v_y,y^5v_y,\\
                                                                                                                                                        &&&y^2u_{xx},yv_{xx},y^4p_{xx},y^3u_{xy},y^2v_{xy},y^5p_{xy},y^4u_{yy},\\
                                                                                                                                                        &&&y^3v_{yy},y^6p_{yy} \\
  \hline\end{array}\end{eqnarray*}\end{table}
 where * and ** are refer to ordinary invariants and first order
differential invariants respectively.
\section{Characterization of invariant system of differential equations}
Differential invariants help us to find general systems of
differential equations which admit a prescribed symmetry group.
One say, if $G$ is a symmetry group for a system of PDEs with
functionally differential invariants, then, the system can be
rewritten in terms of differential invariants. For finding the
differential invariants of the system (\ref{eq:22}) up to order 2,
we should solve the following systems of PDEs:
\begin{eqnarray}\label{eq:11}
\begin{array}{lclclclcl}
\displaystyle\frac{\partial I}{\partial
x}=0,&&\displaystyle\frac{\partial I}{\partial
y}=0,&&\displaystyle\frac{\partial I}{\partial
p}=0,&&\displaystyle x\frac{\partial I}{\partial
x}+u\frac{\partial I}{\partial u}+2p\frac{\partial I}{\partial
p}=0,&&\displaystyle y\frac{\partial I}{\partial
y}-2u\frac{\partial I}{\partial u}-v\frac{\partial I}{\partial
v}-4p\frac{\partial I}{\partial p}=0,
\end{array}
\end{eqnarray}
where $I$ is a smooth function of $(x,y,u,v,p)$,
\begin{eqnarray}\label{eq:12}
\begin{array}{lclcl}
&&\displaystyle\frac{\partial I_1}{\partial
x}=0,\qquad\displaystyle\frac{\partial I_1}{\partial
y}=0,\qquad\displaystyle\frac{\partial I_1}{\partial p}=0,\\
&&\displaystyle x\frac{\partial I_1}{\partial x}+u\frac{\partial
I_1}{\partial u}+2p\frac{\partial I_1}{\partial
p}-v_x\frac{\partial I_1}{\partial
v_x}+\cdots+2p_y\frac{\partial I_1}{\partial p_y}=0,\\
&&\displaystyle y\frac{\partial I_1}{\partial y}-2u\frac{\partial
I_1}{\partial u}-v\frac{\partial I_1}{\partial v}-4p\frac{\partial
I_1}{\partial p}-2u_x\frac{\partial I_1}{\partial
u_x}-\cdots-5p_y\frac{\partial I_1}{\partial p_y}=0,
\end{array}
\end{eqnarray}
where $I_1$ is a smooth function of
$(x,y,u,v,p,u_x,v_x,p_x,u_y,v_y,p_y)$,
\begin{eqnarray}\label{eq:24}
\begin{array}{lclcl}
&&\displaystyle\frac{\partial I_2}{\partial
x}=0,\qquad\displaystyle\frac{\partial I_2}{\partial
y}=0,\qquad\displaystyle\frac{\partial I_2}{\partial p}=0,\\
&&\displaystyle x\frac{\partial I_2}{\partial x}+u\frac{\partial
I_2}{\partial u}+2p\frac{\partial I_2}{\partial
p}-v_x\frac{\partial I_2}{\partial
v_x}+\cdots+2p_{yy}\frac{\partial I_2}{\partial p_{yy}}=0,\\
&&\displaystyle y\frac{\partial I_2}{\partial y}-2u\frac{\partial
I_2}{\partial u}-v\frac{\partial I_2}{\partial v}-4p\frac{\partial
I_2}{\partial p}-2u_x\frac{\partial I_2}{\partial
u_x}-\cdots-6p_{yy}\frac{\partial I_2}{\partial p_{yy}}=0.
\end{array}
\end{eqnarray}
The solutions of (\ref{eq:11}), (\ref{eq:12}) and (\ref{eq:24})
are
\begin{eqnarray}
I&=&const.,\nonumber\\
\displaystyle
I_1&=&\Phi\Big(\frac{u_x}{v^2},\frac{uv_x}{v^3},\frac{p_x}{uv^2},\frac{u_y}{uv},\frac{v_y}{v^2},\frac{p_y}{u^2v}\Big),\label{eq:14}\\
\displaystyle
I_2&=&\Psi\Big(\frac{u_x}{v^2},\frac{uv_x}{v^3},\frac{p_x}{uv^2},\frac{u_y}{uv},\frac{v_y}{v^2},\frac{p_y}{u^2v},\frac{uu_{xx}}{v^4},\frac{u^2v_{xx}}{v^5},\frac{p_{xx}}{v_4},\frac{u_{xy}}{v^3},\frac{uv_{xy}}{v^4},\frac{p_{xy}}{uv^3},\frac{u_{yy}}{uv^2},\frac{v_{yy}}{v^3},\frac{p_{yy}}{u^2v^2}\Big),\label{eq:15}
\end{eqnarray}
thus, the most general first and second order system (\ref{eq:22})
admitting (\ref{eq:13}) as a symmetry group can be rewritten with
two arbitrary functions (\ref{eq:14}) and (\ref{eq:15}).
Characterization of invariant system of differential equations is
an alternative method for reduction of differential equations
invariant under a symmetry group.
\section{Optimal system of steady two-dimensional boundary-layer stagnation-point flow equations}
Let a system of differential equation $\Delta$ admitting the
symmetry Lie group $G$,be given. Now $G$ operates on the set of
solutions $S$ of $\Delta$. Let $s\cdot G$ be the orbit of $s$, and
$H$ be a subgroup of $G$. Invariant $H-$solutions $s\in S$ are
characterized by equality $s\cdot S=\{s\}$. If $h\in G$ is a
transformation and $s\in S$,then $h\cdot(s\cdot H)=(h\cdot s)\cdot
(hHh^{-1})$. Consequently,every invariant $H-$solution $s$
transforms into an invariant $hHh^{-1}-$solution (Proposition 3.6
of \cite{[5]}).

$~~~$Therefore, different invariant solutions are found from
similar subgroups of $G$. Thus, classification of invariant
$H-$solutions is reduced to the problem of classification of
subgroups of $G$,up to similarity. An optimal system of
$s-$dimensional subgroups of $G$ is a list of conjugacy
inequivalent $s-$dimensional subgroups of $G$ with the property
that any other subgroup is conjugate to precisely one subgroup in
the list. Similarly, a list of $s-$dimensional subalgebras forms
an optimal system if every $s-$dimensional subalgebra of $\goth g$
is equivalent to a unique member of the list under some element of
the adjoint representation: $\tilde{\goth h}={\rm
Ad}(g)\cdot{\goth h}$.

$~~~$Let $H$ and $\tilde{H}$ be connected, $s-$dimensional Lie
subgroups of the Lie group $G$ with corresponding Lie subalgebras
${\goth h}$ and $\tilde{\goth h}$ of the Lie algebra ${\goth g}$
of $G$. Then $\tilde{H}=gHg^{-1}$ are conjugate subgroups if and
only $\tilde{\goth h}={\rm Ad}(g)\cdot{\goth h}$ are conjugate
subalgebras (Proposition 3.7 of \cite{[5]}). Thus,the problem of
finding an optimal system of subgroups is equivalent to that of
finding an optimal system of subalgebras, and so we concentrate on
it.
\subsection{One-dimensional optimal system}
For one-dimensional subalgebras, the classification problem is
essentially the same as the problem of classifying the orbits of
the adjoint representation, since each one-dimensional subalgebra
is determined by a nonzero vector in Lie algebra symmetries of
steady two-dimensional boundary-layer stagnation-point flow
equations and so to "simplify" it as much as possible.

The adjoint action is given by the Lie series
\begin{eqnarray}\label{eq:9}
\mbox{Ad}(\exp(\varepsilon\textbf{v}_i)\textbf{v}_j)=\textbf{v}_j-\varepsilon[\textbf{v}_i,\textbf{v}_j]
+\frac{\varepsilon^2}{2}[\textbf{v}_i,[\textbf{v}_i,\textbf{v}_j]]-\cdots,
\end{eqnarray}
where $[\textbf{v}_i,\textbf{v}_j]$ is the commutator for the Lie
algebra, $\varepsilon$ is a parameter, and $i,j=1,\cdots,5$. Let
$F^{\varepsilon}_i:{\goth g}\rightarrow{\goth g}$ defined by
$\textbf{v}\mapsto\mbox{Ad}(\exp(\varepsilon\textbf{v}_i)\textbf{v})$
is a linear map, for $i=1,\cdots,5$. The matrices
$M^\varepsilon_i$ of $F^\varepsilon_i$, $i=1,\cdots,5$, with
respect to basis $\{\textbf{v}_1,\cdots,\textbf{v}_5\}$ are
\begin{eqnarray}\label{eq:6}
\begin{array}{lclcl}
M^\varepsilon_1=\small\left(\begin{array}{ccccc}
1&0&0&0&0\\0&1&0&0&0\\0&0&1&0&0\\-\varepsilon&0&0&1&0\\0&0&0&0&1\end{array}
\right),&& M^\varepsilon_2=\small\left(\begin{array}{ccccc}
1&0&0&0&0\\0&1&0&0&0\\0&0&1&0&0\\0&0&0&1&0\\0&-\varepsilon&0&0&1\end{array}
\right),&& M^\varepsilon_3=\small\left(\begin{array}{ccccc}
1&0&0&0&0\\0&1&0&0&0\\0&0&1&0&0\\0&0&-2\varepsilon&1&0\\0&0&4\varepsilon&0&1\end{array}
\right),\\ M^\varepsilon_4=\small\left(\begin{array}{ccccc}
e^\varepsilon&0&0&0&0\\0&1&0&0&0\\0&0&e^{2\varepsilon}&1&0\\0&0&0&1&0\\0&0&0&0&1\end{array}
\right),&& M^\varepsilon_5=\small\left(\begin{array}{ccccc}
1&0&0&0&0\\0&e^\varepsilon&0&0&0\\0&0&e^{-4\varepsilon}&0&0\\0&0&0&1&0\\0&0&0&0&1\end{array}
\right)
\end{array}
\end{eqnarray}
by acting these matrices on a vector field $\textbf{v}$
alternatively we can  show that a one-dimensional optimal system
of ${\goth g}$ is given by
\begin{eqnarray}\label{eq:20}
\begin{array}{lclclcl}
X_1=\textbf{v}_3,&&X_2=\alpha_1\textbf{v}_1+\alpha_2\textbf{v}_2,&&X_3=\alpha_1\textbf{v}_2+\alpha_2\textbf{v}_3,&&X_4=\alpha_1\textbf{v}_1+\alpha_2\textbf{v}_2+\alpha_3\textbf{v}_3,
\end{array}
\end{eqnarray}
where $\alpha_i$'s are real constants.
\subsection{Two-dimensional optimal system}
Next step is to construct two-dimensional optimal system, i.e.,
classification of two-dimensional subalgebras of $\goth g$. The
process is by selecting one of the vector fields in (\ref{eq:20}),
say, any vector field of (\ref{eq:20}). Let us consider $X_1$ (or
$X_i,i=2,3,4,5$). Corresponding to it, a vector field
$X=a_1\textbf{v}_1+\cdots+a_5\textbf{v}_5$, where $a_i$'s are
smooth functions of $(x,y,u,v,p)$ is chosen, so we must have
\begin{eqnarray}\label{eq:4}
[X_1,X]=\lambda X_1+\mu X,
\end{eqnarray}
the equation (\ref{eq:4}) leads us to the system
\begin{eqnarray}\label{eq:5}
C^i_{jk}\alpha_ja_k=\lambda
a_i+\mu\alpha_i\hspace{2cm}(i=1,2,3,4,5).
\end{eqnarray}
The solutions of the system (\ref{eq:5}), give one of the
two-dimensional generator and the second generator is $X_1$ or,
$X_i,i=2,3,4,5$ if selected. After the construction of all
two-dimensional subalgebras, for every vector fields of
(\ref{eq:20}), they need to be simplified by the action of
(\ref{eq:6}) in the manner analogous to the way of one-dimensional
optimal system.

Consequently the two-dimensional optimal system of $\goth g$ has
three classes of $\goth g$'s members combinations such as
\begin{eqnarray}\label{eq:7}
\begin{array}{lclclcl}
\textbf{v}_1,\textbf{v}_2,&&\textbf{v}_1,\textbf{v}_3,&&\textbf{v}_2,\textbf{v}_3,&&\textbf{v}_3,\textbf{v}_4,\\
\textbf{v}_1,\textbf{v}_4,&&\textbf{v}_1,\textbf{v}_5,&&\textbf{v}_2,\textbf{v}_4&&\textbf{v}_2,\textbf{v}_5,\\
\textbf{v}_3,\textbf{v}_5,&&\textbf{v}_4,\textbf{v}_5,&&\textbf{v}_1,\sum_{i=2}^5\beta_i\textbf{v}_i,&&\beta_1\textbf{v}_1+\beta_2\textbf{v}_2,\textbf{v}_3+\beta_3(\textbf{v}_4+\textbf{v}_5),\\
\beta_1\textbf{v}_2+\beta_2\textbf{v}_3,\textbf{v}_1+\frac{5}{2}\beta_3(\textbf{v}_4+\textbf{v}_5).
\end{array}
\end{eqnarray}
\subsection{Three-dimensional optimal system}
This system can be developed by the method of expansion of
two-dimensional optimal system. For this take any two-dimensional
subalgebras of (\ref{eq:7}), let us consider the first two vector
fields of (\ref{eq:7}), and call them $Y_1$ and $Y_2$, thus, we
have a subalgebra with basis $\{Y_1,Y_2\}$, find a vector field
$Y=a_1\textbf{v}_1+\cdots+a_5\textbf{v}_5$, where $a_i$'s are
smooth functions of $(x,y,\bar{u},\bar{v},\bar{p})$, such the
triple $\{Y_1,Y_2,Y\}$ generates a basis of a three-dimensional
algebra. For that it is necessary an sufficient that the vector
field $Y$ satisfies the equations
\begin{eqnarray}\label{eq:8}
[Y_1,Y]&=&\lambda_1Y+\mu_1Y_1+\nu_1Y_2,\nonumber\\
{[Y_2,Y]}&=&\lambda_2Y+\mu_2Y_1+\nu_2Y_2,
\end{eqnarray}
and following from (\ref{eq:8}), we obtain the system
\begin{eqnarray}\label{eq:10}
C^i_{jk}\beta_r^ja_k&=&\lambda_1a_i+\mu_1\beta_r^i+\nu_1\beta_s^i,\hspace{1cm}\alpha=1,2,3,4,5\\
C^i_{jk}\beta_s^ja_k&=&\lambda_2a_i+\mu_2\beta_r^i+\nu_2\beta_s^i,\hspace{1cm}\alpha=1,2,3,4,5.\nonumber
\end{eqnarray}
The solutions of system (\ref{eq:10}) is linearly independent of
$\{Y_1,Y_2\}$ and give a three-dimensional subalgebra. This
process is used for the another two couple vector fields of
(\ref{eq:7}).

Consequently the three-dimensional optimal system of $\goth g$ is
given by
\begin{eqnarray}\label{eq:16}
\begin{array}{lclclclcl}
\textbf{v}_1,\textbf{v}_2,\textbf{v}_3,&&\textbf{v}_1,\textbf{v}_2,\textbf{v}_4,&&\textbf{v}_1,\textbf{v}_2,\textbf{v}_5,&&\textbf{v}_1,\textbf{v}_3,\textbf{v}_4,&&\textbf{v}_1,\textbf{v}_3,\textbf{v}_5,\\
\textbf{v}_1,\textbf{v}_4,\textbf{v}_5,&&\textbf{v}_2,\textbf{v}_3,\textbf{v}_4,&&\textbf{v}_2,\textbf{v}_3,\textbf{v}_5,&&\textbf{v}_2,\textbf{v}_4,\textbf{v}_5,&&\textbf{v}_3,\textbf{v}_4,\textbf{v}_5.
\end{array}
\end{eqnarray}
\subsection{Four-dimensional optimal system}
Four-dimensional optimal system can be developed by the method of
expansion of three-dimensional optimal system analogous to
three-dimensional. First step is selecting a three-dimensional
subalgebra of (\ref{eq:16}), and call them $Z_1,Z_2$ and $Z_3$,
then find a vector field
$Z=a_1\textbf{v}_1+\cdots+a_5\textbf{v}_5$, where $a_i$'s are
smooth functions of $(x,y,\bar{u},\bar{v},\bar{p})$, such the
foursome $\{Z_1,Z_2,Z_3,Z\}$ generates a basis of a
four-dimensional algebra. It is necessary an sufficient that the
vector field $Z$ satisfies the equations
\begin{eqnarray}\label{eq:17}
[Z_1,Z]  &=&\lambda_1Z+\mu_1Z_1+\nu_1Z_2+\gamma_1Z_3,\nonumber\\
{[Z_1,Z]}&=&\lambda_2Z+\mu_2Z_1+\nu_2Z_2+\gamma_2Z_3,\\
{[Z_1,Z]}&=&\lambda_3Z+\mu_3Z_1+\nu_3Z_2+\gamma_3Z_3,\nonumber
\end{eqnarray}
following from (\ref{eq:17}), we obtain the system
\begin{eqnarray}\label{eq:19}
C^i_{jk}a_k&=&\lambda_1a_i+\mu_1+\nu_1+\gamma_1,\hspace{1cm}\alpha=1,2,3,4,5,\nonumber\\
C^i_{jk}a_k&=&\lambda_2a_i+\mu_2+\nu_2+\gamma_2,\hspace{1cm}\alpha=1,2,3,4,5,\\
C^i_{jk}a_k&=&\lambda_3a_i+\mu_3+\nu_3+\gamma_3,\hspace{1cm}\alpha=1,2,3,4,5.\nonumber
\end{eqnarray}
The solutions of system (\ref{eq:19}) is linearly independent of
$\{Z_1,Z_2,Z_3\}$ and give a four-dimensional subalgebra. This
process is used for the another triple vector fields of
(\ref{eq:16}).
All previous calculations lead to the table (\ref{table:4}) for
the optimal system of $\goth g$.
\begin{table}
\caption{Optimal system of subalgebras}\label{table:4}
\vspace{-0.3cm}\begin{eqnarray*}\hspace{-0.75cm}\begin{array}{lclclclclclclclclclclclclclc}
\hline
  \mbox{dimension}   &\hspace{-0.1cm}1                                                                              &\hspace{1cm}2                                                                                                                     &\hspace{0.1cm}3                                                          &\hspace{1.5cm}4                                                                    &\hspace{0.5cm}5\\ \hline
                     &\hspace{-0.1cm} \langle\textbf{v}_3\rangle                                                    &\hspace{0.5cm} \langle\textbf{v}_1,\textbf{v}_2\rangle                                                                            &\hspace{0.1cm}\langle\textbf{v}_1,\textbf{v}_2,\textbf{v}_3\rangle       &\hspace{0.5cm}\langle\textbf{v}_1,\textbf{v}_2,\textbf{v}_3,\textbf{v}_4\rangle    &\hspace{0.5cm}\langle\textbf{v}_1,\textbf{v}_2,\textbf{v}_3,\textbf{v}_4,\textbf{v}_5\rangle\\
  \mbox{subalgebras} &\hspace{-0.1cm} \langle\alpha_1\textbf{v}_1+\alpha_2\textbf{v}_2\rangle                       &\hspace{0.5cm} \langle\textbf{v}_1,\textbf{v}_3\rangle                                                                            &\hspace{0.1cm}\langle\textbf{v}_1,\textbf{v}_2,\textbf{v}_4\rangle       &\hspace{0.5cm}\langle\textbf{v}_1,\textbf{v}_2,\textbf{v}_3,\textbf{v}_5\rangle    &\hspace{0.5cm}                                                                              \\
                     &\hspace{-0.1cm} \langle\alpha_1\textbf{v}_2+\alpha_2\textbf{v}_3\rangle                       &\hspace{0.5cm} \langle\textbf{v}_2,\textbf{v}_3\rangle                                                                            &\hspace{0.1cm}\langle\textbf{v}_1,\textbf{v}_2,\textbf{v}_5\rangle       &\hspace{0.5cm}\langle\textbf{v}_1,\textbf{v}_2,\textbf{v}_4,\textbf{v}_5\rangle    &\hspace{0.5cm}                                                                              \\
                     &\hspace{-0.1cm} \langle\alpha_1\textbf{v}_1+\alpha_2\textbf{v}_2+\alpha_3\textbf{v}_3\rangle  &\hspace{0.5cm} \langle\textbf{v}_2,\textbf{v}_4\rangle                                                                            &\hspace{0.1cm}\langle\textbf{v}_1,\textbf{v}_3,\textbf{v}_4\rangle       &\hspace{0.5cm}\langle\textbf{v}_1,\textbf{v}_3,\textbf{v}_4,\textbf{v}_5\rangle    &\hspace{0.5cm}                                                                              \\
                     &\hspace{-0.1cm}                                                                               &\hspace{0.5cm} \langle\textbf{v}_3,\textbf{v}_4\rangle                                                                            &\hspace{0.1cm}\langle\textbf{v}_1,\textbf{v}_3,\textbf{v}_5\rangle       &\hspace{0.5cm}\langle\textbf{v}_2,\textbf{v}_3,\textbf{v}_4,\textbf{v}_5\rangle    &\hspace{0.5cm}                                                                              \\
                     &\hspace{-0.1cm}                                                                               &\hspace{0.5cm} \langle\textbf{v}_1,\textbf{v}_4\rangle                                                                            &\hspace{0.1cm}\langle\textbf{v}_1,\textbf{v}_4,\textbf{v}_5\rangle       &\hspace{0.5cm}                                                                     &\hspace{0.5cm}                                                                              \\
                     &\hspace{-0.1cm}                                                                               &\hspace{0.5cm} \langle\textbf{v}_1,\textbf{v}_5\rangle                                                                            &\hspace{0.1cm}\langle\textbf{v}_2,\textbf{v}_3,\textbf{v}_4\rangle       &\hspace{0.5cm}                                                                     &\hspace{0.5cm}                                                                              \\
                     &\hspace{-0.1cm}                                                                               &\hspace{0.5cm} \langle\textbf{v}_2,\textbf{v}_5\rangle                                                                            &\hspace{0.1cm}\langle\textbf{v}_2,\textbf{v}_3,\textbf{v}_5\rangle       &\hspace{0.5cm}                                                                     &\hspace{0.5cm}                                                                              \\
                     &\hspace{-0.1cm}                                                                               &\hspace{0.5cm} \langle\textbf{v}_3,\textbf{v}_5\rangle                                                                            &\hspace{0.1cm}\langle\textbf{v}_2,\textbf{v}_4,\textbf{v}_5\rangle       &\hspace{0.5cm}                                                                     &\hspace{0.5cm}                                                                              \\
                     &\hspace{-0.1cm}                                                                               &\hspace{0.5cm} \langle\textbf{v}_4,\textbf{v}_5\rangle                                                                            &\hspace{0.1cm}\langle\textbf{v}_3,\textbf{v}_4,\textbf{v}_5\rangle       &\hspace{0.5cm}                                                                     &\hspace{0.5cm}                                                                              \\
                     &\hspace{-0.1cm}                                                                               &\hspace{0.5cm} \langle\textbf{v}_1,\sum_{i=2}^5\beta_i\textbf{v}_i\rangle                                                         &\hspace{0.1cm}                                                           &\hspace{0.5cm}                                                                     &\hspace{0.5cm}                                                                              \\
                     &\hspace{-0.1cm}                                                                               &\hspace{0.5cm} \langle\beta_1\textbf{v}_1+\beta_2\textbf{v}_2,\textbf{v}_3+\beta_3(\textbf{v}_4+\textbf{v}_5)\rangle              &\hspace{0.1cm}                                                           &\hspace{0.5cm}                                                                     &\hspace{0.5cm}                                                                              \\
                     &\hspace{-0.1cm}                                                                               &\hspace{0.5cm} \langle\beta_1\textbf{v}_2+\beta_2\textbf{v}_3,\textbf{v}_1+\frac{5}{2}\beta_3(\textbf{v}_4+\textbf{v}_5)\rangle   &\hspace{0.1cm}                                                           &\hspace{0.5cm}                                                                     &\hspace{0.5cm}                                                                              \\
  \hline\end{array}\end{eqnarray*}\end{table}
\section{Lie Algebra Structure}
$\goth g$ has a no any non-trivial \textit{Levi decomposition} in
the form of ${\goth g}={\goth r}\ltimes{\goth g}_1$, because
$\goth g$ has no any non-trivial radical, i.e., if $\goth r$ be
the radical of $\goth g$, then ${\goth g}={\goth r}$.

$~~~$If we want to integration an involuting distribution, the
process decomposes into two steps:
\begin{itemize}
\item integration of the evolutive distribution with symmetry Lie
algebra ${\goth g}/{\goth r}$, and
\item integration on integral manifolds with symmetry algebra $\goth
r$.
\end{itemize}
First, applying this procedure to the radical $\goth r$ we
decompose the integration problem into two parts: the integration
of the distribution with semisimple algebra ${\goth g}/{\goth r}$,
then the integration of the restriction of distribution to the
integral manifold with the solvable symmetry algebra $\goth r$.\\

The last step can be performed by quadratures. Moreover, every
semisimple Lie algebra ${\goth g}/{\goth r}$ is a direct sum of
simple ones which are ideal in ${\goth g}/{\goth r}$. Thus, the
Lie-Bianchi theorem reduces the integration problem to involutive
distributions equipped with simple algebras of symmetries. Thus,
integrating of system (\ref{eq:22}), become so much easy.

$~~~$The Lie algebra $\goth g$ is solvable and non-semisimple. It
is solvable because if ${\goth
g}^{(1)}=\langle\textbf{v}_i,[\textbf{v}_i,\textbf{v}_j]\rangle=[\goth
g, \goth g]$, we have ${\goth g}^{(1)}=[{\goth g},{\goth g}]
=\langle \textbf{v}_1,\cdots,\textbf{v}_5\rangle$, and ${\goth
g}^{(2)}=[{\goth g}^{(1)},{\goth g}^{(1)}] =\langle
\textbf{v}_1,\textbf{v}_2,2\textbf{v}_3\rangle$, so, we have a
chain of ideals ${\goth g}^{(1)}\supset{\goth
g}^{(2)}\supset\{0\}$, and it is non-semisimple because its
killing form
\begin{eqnarray*}
\left(\begin{array}{ccccc}
0&0&0&0&0\\0&0&0&0&0\\0&0&0&0&0\\0&0&0&5&-8\\0&0&0&-8&17\end{array}
\right),
\end{eqnarray*}
is degenerate.

$~~~$According to the table of the commutators, $\goth g$ has two
abelian 2 and 3-dimensional subalgebras spanned by
$\langle\textbf{v}_1,\textbf{v}_2,\textbf{v}_3\rangle$ and
$\langle\textbf{v}_4,\textbf{v}_5\rangle$ respectively such that
the first one is an ideal in $\goth g$, so, we can decompose
$\goth g$ in the to semidirect sum of ${\goth
g}=\Bbb{R}^3\ltimes\Bbb{R}^2$.
\section{Conclusion}
In this article group classification of Newtonian incompressible
fluid's equations flow in turbulent boundary layers and the
algebraic structure of the symmetry group is considered.
Classification of r-dimensional subalgebra is determined by
constructing r-dimensional optimal system. Some invariant objects
are fined and the Lie algebra structure of symmetries is found.


\begin{thebibliography}{9}
\bibitem{[1]}{\sc M.L. Gandarias and M.S. Bruzon}, {Type II didden symmetries through weak symmetries for some wave
equation}, {Communications in Nnonlinear Science and Numerical
Simulation}, 2009, article in press.
\bibitem{[2]}{\sc M. Nadjafikhah and S.R. Hejazi}, {Symmetry analysis of cylindrical Laplace
equation}, {Balkan journal of geometry and applications}, 2009.
\bibitem{[3]}{\sc M. Nadjafikhah, R. Bakhshandeh-Chamazkotia},  and {\sc A.
Mahdipour-Shirayeha}, {A symmetry classification for a class of
(2+1)-nonlinear wave equation}, Article in press,
doi:10.1016/j.na.2009.03.087.
\bibitem{[4]}{\sc P.J. Olver}, {Equivalence, Invariant and Symmetry},
{Cambridge University Press}, {Cambridge University Press,
Cambridge 1995.}
\bibitem{[5]}{\sc P.J. Olver}, {Applications of Lie Groups to
Differential equations}, {Second Edition, GTM, Vol. 107, Springer
Verlage, New York, 1993.}
\bibitem{[6]}{\sc L.V. Ovsiannikov}, {Group Analysis of
Differential Equations}, {Academic Press, New York, 1982.}
\bibitem{[7]} {\sc G.B. Whitham}, {Linear and Nonlinear Waves. New York: Wiley, p. 617, 1974.}
\bibitem{[8]}{\sc D. Zwillinger}, {Handbook of Differential Equations, 3rd ed. Boston, MA: Academic Press, p. 132,
1997.}
\bibitem{[9]}{\sc R. Cimpoiasu and R. Constantinescu}, {Lie symmetries and invariants for a 2D nonlinear heat equation, Nonlinear analysis, 68 (2008) 2261-1168.}
\bibitem{[10]}{\sc Q. Changzheng and Q. Huang}, {Symmetry reduction and exact solution of the affine heat equation, J. Math. Appl. 346 (2008) 521-530.}
\bibitem{[11]}{\sc P.G.L. Leach}, {Symmetry and singularity properties of the generalised Kummer-Schwartz and related equations, J. Math. Appl. 348 (2008) 467-493.}
\bibitem{[12]}{\sc H. Azad, M.T. Mustafa}, {Symmetry analysis of wave equation on sphere, J. Math. Appl. 333 (2007) 1180-1188.}
\bibitem{[13]}{\sc E. Demetriou, N.M. Ivanova, C. Sophocleous}, {Group analysis of (2+1)- and (3+1)-dimensional diffusion-convection equations, J. Math. Appl. 348 (2008) 55-65.}
\bibitem{[14]}{\sc A. Kushner, V. Lychagin and V. Rubstov}, {Contact geometry and nonlinear differential equations, Cambridge University Press,
Cambridge 2007.}
\bibitem{[13]}\verb"www.wikipedia.com".
\end{thebibliography}
\end{document}